\documentclass[11pt]{amsart}\usepackage{amsfonts} \usepackage{amscd}
\usepackage{amssymb}

\newcommand{\ring}[1]{\mathbb{#1}}
\newcommand{\C}{\ring{C}} \newcommand{\Q}{\ring{Q}}
\newcommand{\Z}{\ring{Z}}

\newcommand{\boldL}{{\bf L}}
\newcommand{\lef}{\ring{L}}

\newcommand{\be}{\begin{equation}}
\newcommand{\ee}{\end{equation}}
\newcommand{\nd}{\noindent}

\def\1{{\mu\mkern-6mu\mu}}

\def\op#1{{\operatorname{#1}}}
\newcommand{\lie}{{\bf\mathfrak g}}
\newcommand{\lieh}{{\bf\mathfrak h}}

\newcommand{\gl}{{\bf\mathfrak {gl}}}
\newcommand{\so}{{\bf\mathfrak {so}}}
\newcommand{\syp}{{\bf\mathfrak {sp}}}
\newcommand{\un}{{\bf\mathfrak u}}
\newcommand{\orb}{{\mathcal O}}
\newcommand{\exi}{{\exists}}
\newcommand{\sign}{{\bf\text{sign}}}
\newcommand{\sgn}{{\text{sgn}}}

\newcommand{\lin}{{\text{lin}}}
\newcommand{\res}{{\text{res}}}
\newcommand{\li}{{\text{lin.ind}}}
\newcommand{\monic}{{\text{monic}}}
\newcommand{\irred}{{\text{irred}}}
\newcommand{\Frob}{{\text{Frob}}}
\def\rq{\text{'}}

\title{Virtual Transfer Factors}
\author{Julia Gordon and Thomas C. Hales}

\date{August 23, 2002}

\begin{document}

\theoremstyle{plain}
\newtheorem{thm}{Theorem}
\newtheorem{lem}[thm]{Lemma}
\newtheorem{cor}[thm]{Corollary}
\newtheorem{prop}[thm]{Proposition}

\theoremstyle{definition}
\newtheorem{rem}[thm]{Remark}
\newtheorem{defn}[thm]{Definition}
\newtheorem{ex}[thm]{Example}

\begin{abstract}  The Langlands-Shelstad transfer factor is a
function defined on some reductive groups over a $p$-adic field.
Near the origin of the group, it may be viewed as a function on
the Lie algebra.  For classical groups, its values have the form
$q^c\,\sign$, where $\sign\in\{-1,0,1\}$, $q$ is the cardinality
of the residue field, and $c$ is a rational number. The $\sign$
function partitions the Lie algebra into three subsets. This
article shows that this partition into three subsets is
independent of the $p$-adic field in the following sense.  We
define three universal objects (virtual sets in the sense of
Quine) such that for any $p$-adic field $F$ of sufficiently large
residue characteristic, the $F$-points of these three virtual sets
form the partition.

The theory of arithmetic motivic integration associates a virtual
Chow motive with each of the three virtual sets.  The construction
in this article achieves the first step in a long program to
determine the (still conjectural) virtual Chow motives that
control the behavior of orbital integrals.
\end{abstract}
\maketitle

\setcounter{section}{-1}
\section{Introduction}

\subsection{The Langlands-Shelstad transfer factor}

Langlands and Shelstad have introduced a function, called the {\it
transfer factor}, on certain reductive groups over a $p$-adic
field. The definition of the transfer factor involves the theory
of endoscopy and various constructs of local class field theory.

The transfer factor is expected to figure prominently in the
development of the Langlands program.  The {\it fundamental
lemma}, a conjectural system of identities between orbital
integrals, is expressed by means of transfer factors.  A proof of
the fundamental lemma is needed for many applications of the trace
formula.

The definition of the transfer factor has been simplified in
special contexts by Hales \cite{H}, Kottwitz \cite{K}, and
Waldspurger \cite{W}.  This article follows Waldspurger's
treatment of transfer factors and draws heavily from \cite{W}.
 He gives a simple definition of
transfer factors in the special case of classical groups defined
by quadratic or hermitian forms.  He makes mild restrictions on
the characteristic of the residue field.   Near the identity
element of the group, the transfer factor can be expressed as a
function on the Lie algebra.  Waldspurger gives an elementary
definition of the transfer factor, as a function on the Lie
algebra.  He proves that it is equivalent to the
Langlands-Shelstad definition (\cite[X]{W}). 
 The list of
classical Lie algebras that we consider appears in Section
\ref{liealgebras}.  Waldspurger's definition will be recalled in
Section \ref{thefactor}.

In the case of the transfer factors that we consider, each value
of the transfer factor is of the form
    $$q^c\,\sign,$$
where $\sign\in\{-1,0,1\}$, $q$ is the order of the residue field,
and $c$ is a rational number.  The definition of the rational
number $c$ is elementary, given as the valuation of a certain
explicit discriminant factor.  The entire complexity of the
transfer factor resides in its sign.

The main purpose of this article is to show that the
Langlands-Shelstad transfer factor on the classical Lie algebras
is given by a formula in the first order language of rings.  At
first glance, the result may appear to be obvious. (Can we not
view the  Langlands-Shelstad definition as a formula for the
transfer factor?)  The force of our result comes from the
restrictive nature of the language we use.  Most of the
fundamental structures of $p$-adic analysis cannot be expressed in
this language. In this language, there are no field extensions or
residue fields, no Galois theory or local class field theory, no
functions apart from addition and multiplication, no additive or
multiplicative characters, no valuation or norm, and no
uniformizing elements. In fact, there are no sets at all in this
language.  It is remarkable that the Langlands-Shelstad transfer
factor can be expressed without any reference to Zermelo Fraenkel
set theory.

Set theory is so entrenched in our usual way of talking about
harmonic analysis on $p$-adic groups, that we are forced to make a
long series of preliminary statements about the ``set-free''
definition of standard constructs such as Lie algebras,
centralizers, orbits of elements, linear spaces, bases, norms,
projection operators, and so forth.

The theory of arithmetic motivic integration allows us to
associate a Chow motive over $\Q$ to formulas in the first order
language of rings.  In this way, we show that the
Langlands-Shelstad transfer factor is {\it motivic} in the sense
of Section \ref{mot_int}. The construction in this article achieves
the first step in a lengthy program to determine the (still
conjectural) virtual Chow motives that control the behavior of
orbital integrals.  The influence of Denef and Loeser's work on
motivic integration should be apparent throughout this article
\cite{DL}.

\subsection{The first order language of rings}

The first order language of rings is a formal language in the
first order predicate calculus.  The concepts of logic and model
theory that we require in this article can be found in Enderton
\cite{E} or Fried and Jarden, \cite{FJ}.  In brief, each element
of a language is a finite sequence of letters from a fixed
alphabet. The letters of the alphabet include countably many
variable symbols $x_i$, constants symbols $c_k$ indexed by a set
$K$, symbols for equality, negation, disjunction, existential
quantification, comma, parentheses, and brackets:
$$
    \begin{matrix}
    = & \neg & \vee & \exists&\\
    ,&( &) &[ & ]\\
    \end{matrix}
$$

There are additional letters in the alphabet for each of a
specified set of function symbols and relation symbols.   The
finite sequences in the language are called expressions.  For a
finite sequence of letters in the alphabet to belong to the
language, it must satisfy various syntactic constraints. The
syntax is built in stages: symbols combine in terms, terms combine
into atomic formulas, atomic formulas combine into formulas.

In the case of the first order language of rings, there are two
constants ($0$ and $1$) and two binary function symbols $+$
(addition) and $\times$ (multiplication).

We allow familiar abbreviations in writing formulas in the
language of rings.  We drop parentheses when they can be
reinserted unambiguously.   We write $+$ as an infix operator
rather than in the usual prefix notation of first-order logic:
    $$`0+1\rq \text{ for } `+(0,1)\rq.$$
We write $`2\rq$, $`3\rq$, and so forth for $`1+1\rq$,
$`1+1+1\rq$, and so forth. We drop the multiplication symbol and
indicate multiplication as a juxtaposition of terms.  Additive inverses
may be introduced: $`a-b\rq$ for $`a+(-1)b\rq$, where $`-1\rq$ is
given through an existential quantifier
\begin{equation}\label{-1}
\exists x (x+1=0).
\end{equation}
(Every formula with $`-1\rq$ can be rewritten without $`-1\rq$ by
replacing $-1$ with a variable symbol $x$ and conjoining the given
formula with (\ref{-1}).)  Division may be introduced in a similar
manner. We use standard logical abbreviations such as the
universal quantifier, conjunction, implication, biconditional:
    $$
    \begin{matrix}
    `\forall x\rq & \text{for} & `\neg\exists x\neg\rq\\
    `a \wedge b\rq & \text{for} & `\neg(\neg a \vee \neg b)\rq\\
    `\phi\Rightarrow\psi\rq & \text{for} & `\neg\phi\vee\psi\rq\\
    `\phi\Leftrightarrow\psi\rq & \text{for} & `(\phi\Rightarrow\psi)
    \wedge (\psi\Rightarrow\phi)\rq\\
    \end{matrix}
    $$

We often use variable symbols that are more suggestive of meaning
than the variable symbols $x_i$ provided by the language.  For
example, if the context is an $n\times n$ matrix, we use variable
symbols $x_{ij}$, $y_{ij}$, and so forth rather than labeling the
variable symbols sequentially $x_1,x_2,\ldots$.  We will
occasionally use a multi-indexing notation.  For example if $X$ is
the matrix $(x_{ij})$ of variable symbols, then $\exists X$ is an
abbreviation of
    $$\exists x_{11}\exists x_{12}\cdots \exists x_{nn}.$$

We write $\phi(x_1,\ldots,x_n)$ to indicate a formula in the first
order language of rings such that all free variables are among the
variable symbols $x_1,\ldots,x_n$.  We use a multi-index notation
here as well, for instance, writing $\phi(X)$, for
$\phi(x_{11},\ldots, x_{nn})$, when $X$ is a matrix of variable
symbols $x_{ij}$.

\subsection{First order language of rings with involution}

We will also have occasion to use the first order language of
rings with involution.  It is constructed in the same way as the
first order language of rings, except that it has an additional
unary function symbol. We write this function symbol as a bar
$\bar t$ over the term $t$ to which the function symbol is
applied.

\subsection{Virtual sets}

The language of first order rings is a highly restrictive language
with no notion of sets.  In particular, the set membership
predicate $\in$ is absent.   Following Quine, we introduce
{\it virtual sets} into the language as abbreviations of various
logical formulas.\footnote{Quine himself calls them {\it virtual
classes}. He tends to use the word `class' in contexts that
mathematicians prefer the word `set.'  In Quine's system,
``Basically, `set' is simply a synonym of `class' that happens to
have more currency than `class' in mathematical contexts.$\ldots$
My own tendency will be to favor the word `class' where `class' or
`set' would do, except for calling the subject set theory
\cite[pages3-4]{Q}.'' Quine's virtual classes and classes are
related to what others call classes and sets,
respectively. See, for example, \cite[page 10]{TZ}.
}
Let $\phi$ be a formula in the first order language of rings.  We
write
    $$`y\in \{x : \phi(x)\}\rq \text{ for } `\phi(y).\rq$$
The construct $\{x: \phi(x)\}$ is called a {\it virtual set}.
Here, $x$ is allowed to be a multi-variable symbol:
    $x = (x_1,\ldots,x_n)$, so that we have
    $$`(y_1,\ldots,y_n)\in \{(x_1,\ldots,x_n) : \phi(x_1,\ldots,x_n)\}\rq\text { for }
    `\phi(y_1,\ldots,y_n)\rq
    $$
When we write $x\in {\bf A}$, it is to be understood that $x$ is a
vector of variable symbols, and that the length of that vector is
the number of free variables in the defining formula of ${\bf A}$.

If $\bf A$ and $\bf B$ are virtual sets defined by formulas $\phi(x)$ and
$\psi(x)$ respectively, we have a notion of subset, union, and
intersection:
    $$
    \begin{matrix}
    `{\bf A}\subset {\bf B}\rq &\text{ for } &`\forall x
    (\phi(x)\Rightarrow\psi(x))\rq\\
    `{\bf A} \cap {\bf B}\rq &\text{ for } &`\{x:\phi(x)\wedge \psi(x)\}\rq\\
    `{\bf A} \cup {\bf B}\rq &\text{ for } &`\{x:\phi(x)\vee\psi(x)\}\rq\\
    \end{matrix}
    $$
In the language there is a single sort of quantifier, a quantifier
of ring sort.  It is impermissible to write an expression such as
$\forall {\bf A}$, where $\bf A$ is a virtual set. (There are no variables
ranging over virtual sets; variables range over elements of
virtual sets.)  The letters `$\bf A$' and `$\bf B$' above are
meta-variables, which lie outside the formal language.
Furthermore, it is impermissible to write one virtual set as an
element of another.

Let $\phi(y,x)$ be a formula with free variables limited to
$y=(y_1,\ldots,y_n)$ and $x = (x_1,\ldots,x_k)$.  We define a {\it
virtual set with parameters\/} $x$ by
    $$`u \in \{y : \phi(y,x)\}\rq \text{ for } `\phi(u,x)\rq$$
where $u = (u_1,\ldots,u_n)$. (The usual cautions about the
capture of free variables apply here.)  We may speak of inclusion,
intersections, and unions of virtual sets with parameters.

In this article, all formulas are taken to be formulas in the
first order language of rings (or rings with involution).  All
virtual sets are understood as given by formulas in this language.

\begin{rem}
It is customary practice to adopt a realist point-of-view in the
discussion of set theory.  That is, mathematical discourse is
framed as a discussion of sets as things, rather than as
well-formed expressions in a formal language.
  We follow a similar practice in this
paper with respect to the well-formed expressions in our formal
language, and adopt a realist stance.  That is, we write this
paper as if the formal language names objects (such as Lie
algebras, centralizers, and orbits).
\end{rem}

\section{Virtual Transfer Factors}

\subsection{Linear algebra}

\begin{defn}  If $\bf V$ is a virtual set, we let
    $\lin({\bf V})$ be the formula
    $$\forall\lambda_1\forall\lambda_2 \, \forall x_1\forall
    x_2\,\,
        (x_1,x_2\in V \Rightarrow \lambda_1 x_1+\lambda_2 x_2\in {\bf V}).$$
That is, $\lin({\bf V})$ asserts that $\bf V$ is a linear space. Here
$\lambda_1$ and $\lambda_2$ are variable symbols and $x_1$ and
$x_2$ are vectors of variable symbols.  The length of the vectors
must equal the number of free variables of the virtual set $\bf V$.
\end{defn}

\begin{defn} \label{independence} (linear independence)
    If $\bf V$ is a virtual set with $N$ free variables, and if
    $e_i$, for $i=1,\ldots,n$, are vectors of terms,
    where each vector has length $N$, we let
    $\li(e_1,\ldots,e_n,{\bf V})$ be the formula
    $$\forall c_1,\ldots,c_n (\sum_{i=1}^n c_i e_i =0\Rightarrow
    c_1=\dots=c_n=0)$$
That is, the formula asserts the linear independence of the
elements $e_1,\ldots,e_n$ in $\bf V$.
\end{defn}

\begin{defn} \label{span} (span)
    If $\bf V$ is a virtual set, we let
    $\op{span}(e_1,\ldots,e_n,{\bf V})$ be the formula
    $$\forall v\in {\bf V} \ \exi
    \lambda_1,\dots,\lambda_n \,\,( v=\sum_{i=1}^n\lambda_ie_i).$$
The formula asserts that $e_1,\ldots,e_n$ span $\bf V$. The length of
the vectors $e_i$ must equal the number of free variables in $\bf V$.
We write $\op{basis}(e_1,\ldots,e_n,{\bf V})$ for the conjunction
    $$\li(e_1,\ldots,e_n,{\bf V}) \wedge \op{span}(e_1,\ldots,e_n,{\bf V}).$$
\end{defn}

\begin{defn} \label{dim n}
If $\bf V$ a virtual set, we let
    $\dim({\bf V},n)$ be the conjunction of the two formulas
    $$
    \begin{array}{lll}
    \lin({\bf V}), \\
    \exi e_1,\dots, e_n\in {\bf V} : \op{basis}(e_1,\ldots,e_n,{\bf V}).
    \end{array}
    $$
That is, $\dim({\bf V},n)$ asserts that $\bf V$ is a linear space of
dimension $n$.
\end{defn}

\begin{defn} \label{Lefschetz}
We let ${\boldL}^n$ be the virtual set
    $$\boldL^n = \{(x_1,\ldots,x_n): \phi(x_1,\ldots,x_n)\}$$
where $\phi$ is any formula with free variables $x_1,\ldots,x_n$
such that $$\forall x_1\cdots \forall x_n (\phi(x_1,\ldots,x_n))$$
is valid. For example, take
    $$`\phi(x_1,\ldots,x_n)\rq\text { to be }
        `(x_1=x_1) \wedge\cdots\wedge (x_n=x_n)\rq.$$
We view this virtual set as the standard linear space of dimension
$n$.   (We use the notation $\boldL^n$ because under the
Denef-Loeser map from formulas to motives, the virtual set
$\boldL^n$ is mapped to the $n$th power of the Lefschetz motive
$\lef$.)
\end{defn}

\subsection{Polynomials}

We distinguish between two types of polynomials.  Polynomials in
the variable symbols appear as terms in the first order language
of rings.  This type of polynomial is not of direct interest to us
in this subsection.

The second type of polynomial is a polynomial in a meta-variable
with coefficients that are terms in the first order language of
rings.  The properties of such polynomials are developed in this
subsection.

We will have need to work with polynomials whose coefficients are
terms in the first order language of rings. We let the $n+1$-tuple
of terms, $(a_0,\ldots,a_n)$, represent the polynomial
    $$\sum a_i \lambda^i.$$
Here, $\lambda$ is a meta-variable, serving as a place holder for
the terms $a_i$.   We may then identify the virtual set
$\monic(n)$ of monic polynomials of degree $n$ with $\boldL^n$.
There is no virtual set of polynomials of arbitrary degree.  If
$f$, $f_1$, and $f_2$ are monic polynomials, then we write
    $$f = f_1 f_2$$
for the conjunction of identities obtained by equating
coefficients:
    $$a_k = \sum_{k=i+j} b_i c_j.$$
Similarly, we write
    $$f = f_1 f_2 \cdots f_\ell$$
for a conjunction of identities of coefficients.

An expression such as
    $$\exists f_1 \exists f_2\quad f = f_1 f_2,$$
for monic polynomials $f$, $f_1$, $f_2$ is to be interpreted as a
disjunction
    $$\bigvee_{n_1+n_2=n} \exists f_1\in\monic(n_1)\quad \exists f_2\in\monic(n_2)
        \quad (f = f_1 f_2).$$
A constraint on the degrees such as
    $$f = f_1 f_2 \text{ with } \deg(f_1)\ge 2$$
should be interpreted as a constraint on the corresponding
disjunction:
    $$\bigvee_{n_1+n_2=n,\,\,n_1\ge2}$$
The formula
    $$\exists f\in\monic(n)\quad \phi(f)$$
is itself to be interpreted as a statement about its coefficients
$a_i$:
    $$\exists a_0,\ldots,a_{n-1} \quad\phi(a_0,\ldots,a_{n-1}).$$

\begin{defn} \label{irreducible}
If $f\in\monic(n)$, we let $\irred(f)$ be the negation of the
formula
    $$\exists f_1 \exists f_2,\quad f = f_1 f_2\text{ with } \deg f_1\ge1,\quad \deg f_2\ge1.$$
$\irred(f)$ asserts the irreducibility of $f$. Note that for each
$n$, $\irred(f)$ is a different formula.  In particular the number
of free variables depends on $n$.
\end{defn}

\begin{defn}\label{even} (even)
If
$f\in\monic(2n)$, we let
    $$\op{even-poly}(f)$$ be the formula
        $a_1 = a_3 = \cdots = a_{2n-1}=0$
asserting that $f$ is an even
    polynomial.
\end{defn}

\subsection{Lie algebras}\label{liealgebras}

\begin{defn} \label{gln}
{\bf The virtual Lie algebra $\gl(n)$:} let
    $$\gl(n) = \{(x_{ij}) : \phi(x_{ij})\},$$
where $\phi$ is a formula in $n^2$ free variables such that
    $$\forall x_{ij} (\phi(x_{ij}))$$
is valid, and $i,j$ range from $1$ to $n$.
\end{defn}

In order to define the virtual classical Lie algebras, we consider
a vector of variable symbols (to be understood as the underlying
linear space) and a bilinear form (presented as a matrix of terms
in the language).
 The classical Lie
algebras will be defined as appropriate virtual subsets of the
virtual set of endomorphisms of the linear space.

\begin{defn}\label{Jeven}
 {\bf The virtual Lie algebra $\so(n)$:} If $n$ is odd, let $J=(q_{ij})_{1\le
i,j\le n}$ be given by
    $$q_{ij}=
        \begin{cases}
            (-1)^{i+1}/2, & \text{if } i+j=n+1,\  i\ne j,\\
            0, & \text{if } i+j\ne n+1,\\
            (-1)^{i+1}, & \text{if } i=j=(n+1)/2.\\
        \end{cases}
    $$
If $n$ is even,  set
    $$q_{ij}=
        \begin{cases}
            0, & \text{if } i+j\ne n+1,\\
            (-1)^{i+1}/2, & \text{if } i+j=n+1,\  i<j,\\
            (-1)^{j+1}/2, & \text{if } i+j=n+1,\ j<i.\\
        \end{cases}
    $$
 The
matrix $J$ is to be understood as a matrix of constant symbols in
the formal language.   Let $X$ be an $n\times n$ matrix of
variable symbols $x_{ij}$.  Define the virtual orthogonal Lie
algebra to be the virtual set
$$\so(n)=\{X : {}^t X J + J X = 0\}.$$
\end{defn}

\begin{defn}\label{sp}
{\bf The virtual Lie algebra $\syp(2r)$} is defined in the same
way:
    $$\syp(2r) = \{X : {}^t X J + J X = 0\}$$
by means of the matrix $J=(q_{ij})_{1\le i,j\le r}$ with
$q_{ij}=(-1)^i$ if $i+j=2r+1$, and zero otherwise.
\end{defn}

\begin{defn}\label{unitary}
{\bf The virtual Lie algebra $\un(r)$} is defined similarly. It is
a virtual set in the first order language of rings with involution
$t\mapsto \bar{t}$. Let
    $$\un(r) = \{X: {}^t\bar X J + J X = 0\}$$
where $J$ is the $r\times r$ matrix of constant symbols given by
$J = (q_{ij})$, with $q_{ij} = 2(-1)^{i+1}$, if $i+j=r+1$ and
$q_{ij}=0$, otherwise.
\end{defn}

\begin{rem} In general, we follow Waldspurger closely in our
definitions, including our choices of bilinear forms
\cite[X.3]{W}. However, in the hermitian case, Waldspurger
introduces a $p$-adic element $\eta$ in a quadratic extension of
the $p$-adic field. Such an element is not definable in the
language of rings with involution, and we are forced to make some
slight adjustments in definitions.
\end{rem}

The symbol $\lie$ will be reserved for one of the classical
virtual Lie algebras that we have defined, or for a direct sum of
such algebras. In general, the bold script in notation
indicates that we are talking about virtual sets.

\subsection{Centralizers}
Let $\lie$ be a virtual Lie algebra, and let $X=(x_{ij})\in\lie$.
Let  $P_X$ be the characteristic polynomial
$$P_X(\lambda)=\det(\lambda\, Id-X),$$
where the determinant is expanded as an explicit polynomial in $x_{ij}$.

\begin{defn}\label{rss} (regular semisimple)
The virtual set $\lie^{reg}$ of regular semisimple elements inside each
classical Lie algebra except for the even orthogonal algebra
is the virtual set
defined by the property that the matrix $X$ has distinct
eigenvalues. Explicitly, it is given by the polynomial condition
$\res(P_X, P_X')\neq 0$, where
$\res(f,g)$ stands for the resultant of
two polynomials $f$ and $g$ (see, e.g., \cite{vdw}, Section 5.4).

The virtual set of regular semisimple elements inside $\so(2r)$ is
the union of two virtual subsets:  the set of matrices in
$\so(2r)$ with distinct eigenvalues and the set $\bf B$ of matrices in
$\so(2r)$ with $0$ as an eigenvalue of multiplicity $2$ and other
eigenvalues distinct:
    $${\bf B} = \{X \in \so(2r) : P_X = \lambda^2 f \wedge
    \res(f,f')\ne0 \wedge f(0)\ne 0\}.$$
\end{defn}

\begin{defn}\label{chp}
We may define stable orbits ${\bf O}^{st}(X)$ as virtual
sets with parameters $X\in\lie^{reg}$ as follows.  For the
unitary, symplectic, and odd orthogonal Lie algebras, we define
    $${\bf O}^{st}(X) = \{Y\in \lie : P_X = P_Y\}.$$

In the case of the even orthogonal Lie algebra, if $X \in
\so(2r)^{reg}$, then $J X$ is a skew symmetric matrix of variable
symbols.  A skew matrix has a Pfaffian $\op{pf}(J X)$ \cite[page
627]{G}. We then have
    $${\bf O}^{st}(X) = \{Y\in\lie : P_X = P_Y \wedge \op{pf}(J
    X)=\op{pf}(J Y)\}.$$
\end{defn}

\begin{defn}\label{central}
Define the centralizer ${\bf C}(X)$ depending on the parameter
$X\in\gl(n)$ by
    $$\{ Y\in\gl(n) : X Y -Y X = 0\}.$$
Similarly, for $X\in\lie$, define
    $${\bf C}_\lie(X) = \{ Y\in\lie : X Y - Y X = 0\}.$$
\end{defn}

\subsection{Lie algebras considered}
\label{considered}

Let $\lie\oplus\lieh$ be one of the following virtual Lie
algebras.
    $$
    \begin{array}{lll}
    \so(2r+1) \oplus \so(2a+1) \oplus \so(2b+1), \text{ with } a+b=r,\\
    \syp(2r)\oplus \syp (2a)\oplus \so(2b), \text{ with } a+b = r,\quad (b\ne1)\\
    \so(2r)\oplus \so(2a)\oplus \so(2b),\text{ with } a + b = r,\quad (a\ne1,b \ne1, r\ne1)\\
    \un(n)\oplus \un(a)\oplus \un(b),\text{ with } a+b = n.\\
    \end{array}
    $$
Each Lie algebra is split except in the unitary case. In each
case, the Lie algebra is a sum of three factors.  We write
$(X,Y,Z)$ for an ordered triple of matrices of variable symbols
corresponding to this direct sum decomposition.  We refer to these
four cases as the odd orthogonal, symplectic, even orthogonal, and
unitary cases, respectively.  We write $\lie$ for the first factor
($\so(2r+1)$, $\syp(2r)$, and so forth), and $\lieh$ for the sum
of the last two factors.

\begin{rem}
The origin of this list of Lie algebras is the following.  Let $F$
be a $p$-adic field.  Let $G$ be a classical quasi-split adjoint
group over $F$ and let $H$ be an elliptic endoscopic group of $G$.
Then the Lie algebras listed above are the Lie algebras of
products $G\times H$.  (In general, $H$ is a product of two
factors.)  The list is not exhaustive.  In particular, it
does not include the non-split even orthogonal groups.
\end{rem}

\begin{rem}  As the introduction to this article explains, the
Langlands-Shelstad transfer factor for classical groups is a
function on the (Lie algebra of) $G\times H$ taking values in the
set
    $$\{-1,0,1\} q^{\Q}.$$
The transfer factor thus partitions the Lie algebra of $G\times H$
into three subsets, corresponding to the possible values ($-1$,
$0$, $1$) of the sign.  By definition, the transfer factor is zero
on elements that are not regular semisimple.   We will realize
these three subsets as virtual subsets of the virtual Lie
algebras. These virtual subsets are what we take as the definition
of the {\it virtual transfer factor}.
\end{rem}

\subsection{The nonzero part of the transfer factor}

We define a virtual subset of $\lie\oplus\lieh$ corresponding to
the set on which the transfer factor is nonvanishing.  Let
$\lie^{reg}$ be the virtual subset of regular semisimple elements
of $\lie$.

\begin{defn}\label{pmset}
 The virtual $\pm$-set of the virtual transfer factor is defined
to be the virtual subset $(\lie\oplus\lieh)_\pm$ of
$\lie\oplus\lieh$ given by
    $$\{(X,Y,Z)\in\lie^{reg}\oplus\lieh^{reg} : P_X = P_Y P_Z\}$$
in the unitary and symplectic cases, by
    $$\{(X,Y,Z)\in\lie^{reg}\oplus\lieh^{reg} : \lambda P_X = P_Y P_Z\}$$
in the odd orthogonal case, and by
    $$\{(X,Y,Z)\in\lie^{reg}\oplus\lieh^{reg} :
    P_X = P_Y P_Z \wedge \op{pf}(JX) =
    (-1)^{a b}\op{pf}(JY)\op{pf}(JZ)\}$$
in the even orthogonal case ($\lieh = \so(2a)\oplus\so(2b)$).
The matrices $J$ are of appropriate
size (symmetric, skew, or hermitian, as appropriate), adapted to
the size of the matrices $X,Y,Z$ as given in Section
\ref{liealgebras}.
\end{defn}

The virtual $0$-set $(\lie\oplus\lieh)_0$ is defined to be the
complement of the $\pm$-set in $\lie\oplus\lieh$.

\subsection{Projection operators}

If $X\in\gl(n)$, then we have the virtual set ${\bf C}(X)$ depending on
parameters $X$.  Let $\op{proj}(X)$ be the virtual set (with parameter
$X$)
    $$
    \op{proj}(X) = \{ P = (p_{ij}) : \forall Y\in {\bf C}(X)
(P Y \in {\bf C}(X) \wedge P P Y =  P Y)\}.$$ That is, it is the
virtual set of matrices $P$ (acting on the same underlying linear
space as $\gl(n)$) such that $P$ acts as a projection operator on
${\bf C}(X)$.

\begin{rem} Return for a moment to the world of set theory, rings,
and modules.  Let $T$ be a linear transformation of $\C^n$ with
distinct eigenvalues, let $P_T\in \C[\lambda]$ be the
characteristic polynomial.  Let $\lambda_i$ be an eigenvalue. Let
$P^{(i)}$ be the characteristic polynomial of $T$ divided by the
factor $(\lambda-\lambda_i)$.  Then
    $$P^{(i)}(\lambda)/P^{(i)}(\lambda_i)$$
is a polynomial, which when evaluated at $T$, defines the
projection operator onto the $\lambda_i$-eigenspace. This
polynomial is uniquely characterized modulo multiples of $P_T$ by
this property.

More generally,
    $$\sum_{i\in S} P^{(i)} (\lambda)/P^{(i)}(\lambda_i)$$
yields the projection operator onto the direct sum of the
eigenspaces of $i\in S$. This polynomial is expressed by means of
the resultant $\res$ in the form
    $$\frac{{\Pi(\lambda,f,\tilde f)}}{\res(f,\tilde f)}$$
for some
    $$\Pi(\lambda,f,\tilde f)\in\Z[\lambda,a_1,\ldots,a_s,b_1,\ldots,b_u],$$
where $$f(t) = t^s + a_{s-1}t^{s-1}+\cdots+a_0
    =\prod_{i\in S} (t-\lambda_i),$$
and $$\tilde f(t) = t^u + b_{u-1}t^{u-1}+\cdots+b_0
    =\prod_{i\not\in S} (t-\lambda_i).$$
The polynomial $\Pi(\lambda,f,\tilde f)$ depends on $T$ only
through the coefficients of $f$ and $\tilde f$.
\end{rem}

Return to the world of virtual sets and formal languages.
If $f\in\monic(s)$ and $\tilde f\in\monic(u)$, and $X =
(x_{ij})$ are variable symbols, we have the matrix
    $$\Pi(X,f,\tilde f)$$
whose coefficients are terms in $x_{ij}$, $a_i$ (the coefficients
of $f$), and $b_i$ (the coefficients of $\tilde f$).

\subsection{The $+1$-set of $\sign$}\label{1-set}

In this section we define the virtual $+1$-set of the transfer
factor.  The definition is slightly different for each case.  In
the unitary case, we must include the involution $t\mapsto\bar
t$.  To increase the uniformity of presentation, we define
    $$t\mapsto\bar t$$
to be the identity map whenever $\lie$ is not the unitary Lie
algebra.

\begin{defn}
If $X$ is a matrix of terms, let $\tau$ be the involution
    $$\tau(X) =(J^{-1}) \,({}^t \bar X) J.$$
The matrix $J$ is that which enters into the definition of $\lie$.
\end{defn}

The lie algebra $\lie$ can be identified with the fixed point set
of $X\mapsto -\tau(X)$.  Let $\gl(n)$ be the linear space containing $\lie$ in
a natural way.

\begin{defn}
For each $k\ge1$, let $\op{even}_k$ be a boolean polynomial in $k$
variables that is true iff an even number of the arguments are
false. For example,
    $$\op{even}_2(b_1,b_2) = (b_1\wedge b_2) \vee
        (\neg b_1\wedge \neg b_2).$$
\end{defn}

In the odd orthogonal case, the characteristic polynomial of a
semisimple element is an odd polynomial, but in the symplectic and
even orthogonal cases, the characteristic polynomial is even.
Again, for the sake of uniform presentation, for any
matrix $Z$ of terms, define $$P_{Z,0}(\lambda) =
\begin{cases}
        P_Z(\lambda)/\lambda, & \lie \text{ odd orthogonal}\\
        P_Z(\lambda), & \text{ otherwise}.
        \end{cases}
        $$

\begin{defn} (norms)  Let $X\in \lie^{reg}$ be a matrix of terms,
and $f$ a monic polynomial.  Let
    $\op{norm}(X,f,U)$
be the formula
    $$\begin{array}{lll}
    \exists \tilde f \exists X_1&\in {\bf C}(X) \,\,\\
    &\left(\,
    P_X = f \tilde f \wedge
    U\in C(X) \wedge
    \Pi(X,f,\tilde f)X_1\tau(X_1)=\Pi(X,f,\tilde f) U
    \right).
    \end{array}$$
It asserts that the `$f$-component' of $U\in C(X)$ is a norm.
\end{defn}

\begin{defn} (trace form)
Let $\op{trace-form}(X,c)$ be the formula
    $$
    \begin{array}{lll}
    \exists e_1,\ldots,e_n
    &\forall x_1,\ldots,x_n,x_1',\ldots,x_n'\\
    &\op{basis}(e_1,\ldots,e_n,C(X))
    \wedge
    c\in C(X)\wedge\\
    &\op{trace}(\tau(\sum x_i e_i)(\sum x_j' e_j)c) =
            {}^t\bar x J x'
\end{array}
    $$
Here, $\op{trace}$ is the matrix trace.  The formula asserts that there is
a basis $e_1,\ldots,e_n$ of ${\bf C}(X)$, for which the trace form on
${\bf C}(X)$ (with constant $c\in C(X)$)
is in agreement with the form
defining $\lie$.
\end{defn}

\begin{defn} (even parity)
Let $\phi(f,\ldots)$ be a formula (or more accurately,
a family of formulas indexed by
the degree of $f$) whose free variables include a monic
polynomial $f$.  Let
    $$\op{even-parity}(f,\phi)$$
    be the formula given by
    $$\begin{array}{lll}
    \bigvee_{\ell}\,
        &\exists f_1,\ldots, f_\ell\,&
         f=f_1\cdots f_\ell \wedge\\
        \qquad\wedge_{i=1}^\ell \deg\,f_i\ge 1&
                    \wedge_{i=1}^\ell\irred(f_i)
                 &\wedge
        \ \op{even}_\ell(\phi(f_1,\ldots),\ldots,\phi(f_\ell,\ldots)).
    \end{array}
    $$
It asserts the even parity of the number of irreducible factors
of $f$ that fail to satisfy $\phi$.
\end{defn}

In the case of symplectic and odd orthogonal groups, we define the
formula $\phi(f,\ldots) = \phi(f,X,c')$ to be
    $$\op{even-poly}(f) \implies \op{norm}(X,f,P_X'(X)c').$$
For even orthogonal, we take $\phi(f,X,c')$ to be
    $$\op{even-poly}(f) \implies
    \exists X'\in {\bf C}(X)\quad
    X X'=1 \wedge \op{norm}(X,f,P_X'(X)c'X').$$

If $X\in\un(n)$, then we work systematically with the
characteristic polynomial of $X\epsilon$ rather than that of $X$,
where $\bar \epsilon = -\epsilon$.  The
characteristic polynomial of $X\epsilon$ has the property
that each of its coefficients is fixed by the involution.
For unitary, we take $\phi(f,X,c',\epsilon)$ to be
    $$\op{norm}(X\epsilon,f,P'_X(X) c').$$

For symplectic and orthogonal groups,
we define the virtual $+1$-set $(\lie\oplus\lieh)_+$ to be the
virtual subset of $(\lie\oplus\lieh)_\pm$ given by $(X,Y,Z)$
satisfying the formula
$$
    \begin{array}{lll}
    \exists c,c'\in {\bf C}(X)&\\
    \\
    \quad c\, c' = 1&\wedge\\
    \quad \tau(c) = \chi c&\wedge\\
    \quad \op{trace-form}(X,c)&\wedge\\
    \\
    \quad \op{even-parity}(P_{Z,0},\phi(\cdot,X,c'))\\
    \end{array}
$$
\medskip
The constant $\chi$ is $\pm1$.  It is $+1$ in each case except for
$\lie = \syp(2r)$ and even unitary.  For $\syp(2r)$
and even unitary, take $\chi=-1$.

In
the unitary case, we replace
    `$\op{even-parity}(P_{Z,0},\phi(\cdot,X,c'))$'
    with
    $$
    `\exists \epsilon\quad
            \epsilon\ne0
            \ \wedge\
        \bar\epsilon = -\epsilon
        \ \wedge\
    \op{even-parity}(P_{Z\epsilon},\phi(\cdot,X,c',\epsilon)).
    \rq
    $$

This completes the definition of the virtual $+1$-set.
We define the virtual $-1$-set $(\lie\oplus\lieh)_-$ by the
complement of the $+1$-set in the $\pm$-set.  This completes our
definition of the virtual transfer factor.

\subsection{Structures}

There are structures (in the sense of model theory) for the first
order theory of rings for every $p$-adic field.  Let $F$ be a
$p$-adic field of characteristic zero.  The domain of the
structure is $F$. The binary operations $+$ and $\times$ become
addition and multiplication in $F$. A structure with domain $F$
attaches a set $\bf A(F)$ to every virtual set $\bf A$ (the set of
$F$-points of $\bf A$, so to speak).

If we take the first order theory of rings with involution, then
we have structures corresponding to separable quadratic extensions
of $p$-adic fields $E/F$.  In the domain $E$, the involution
$t\mapsto \bar t$ becomes the nontrivial automorphism of $E/F$. If
$\bf A$ is a virtual set, and $F$ is a $p$-adic field with
uniquely defined unramified quadratic extension $E$,  we write
${\bf A}(F)$ for the elements of the interpretation of $\bf A$
in the domain $E$ with
involution coming from $E/F$.

\begin{rem}
We write ${\bf A}(F)$ (rather than ${\bf A}(E)$) for the sake
of uniformity of notation, as well as to suggest the analogy with
the $F$-points of a variety (such as the $F$-points of the unitary
group splitting over a quadratic extension $E$).
\end{rem}

\subsection{The Main Theorem}

For the even orthogonal Lie algebra, the set of regular elements
$X$ with eigenvalue $0$ is somewhat exceptional, and will be
excluded from the following theorem.  Such elements are excluded
by restricting to elements $X$ whose image in $\gl(n)$ is regular.

\begin{thm}  Assume $F$ is a $p$-adic field of characteristic zero
of sufficiently large residue field characteristic.  Let
$\lie\oplus\lieh$ be one of the Lie algebras introduced in Section
\ref{liealgebras}.
Let $\Delta(X,Y,Z)$ be
the transfer factor. The set
    $$\{(X,Y,Z)\in\lie(F)\oplus\lieh(F) : \Delta(X,Y,Z)=\sigma \,\wedge X\in \gl(n)^{reg} \}$$
equals
    $$\{(X,Y,Z)\in (\lie\oplus\lieh)_\sigma(F) : X\in \gl(n)^{reg}\}$$
for $\sigma\in\{-1,0,1\}$.
\end{thm}

This theorem will be proved in Section \ref{proof}.

\subsection{Motivic interpretation}\label{mot_int}

By the theory of arithmetic motivic integration, developed by
Denef and Loeser in \cite{DL}, we may associate virtual Chow
motives with virtual sets. (The word {\it virtual\/} is used here
in two different senses.)

Let $\bf A$ be a virtual set defined by a formula $\phi$.  The formula
$\phi$ can be viewed as a formula in Pas's language, which is an
extension of the first order language of rings.  This extension
has additional function symbols $\op{ac}$ and $\op{ord}$, corresponding to
the angular component function and the valuation function in
$p$-adic analysis.  For details, see \cite{Pas}.

A construction of Denef and Loeser \cite[Section 6]{DL} attaches a
virtual Chow motive to formulas in Pas's language.   Orbital
integrals are expected to count points on virtual Chow motives
(see \cite{H1}).  These virtual Chow motives are designed to be
independent of the $p$-adic field.   The construction in this
article achieves the first step in a lengthy program to determine
the virtual Chow motives that control the behavior of orbital
integrals.

Denef and Loeser's construction gives the following corollary. A
ring of virtual Chow motives $K_0^v(\op{Mot}_{\Q,\bar\Q})\otimes\Q$
over $\Q$ is defined in
\cite{DL}.  For any number field $K$, and virtual Chow motive $M$,
we obtain $M_K$ by base change from $\Q$ to $K$.  Let $K_v$ be the
completion of $K$ at the place $v$.  Let $O_v$ be the ring of
integers of $K_v$. The virtual set $(\lie\oplus\lieh)_\pm$ is a
subvariety (in the sense of being defined by polynomial equations)
of $\lie\oplus\lieh$ and it follows that there is a canonically
defined Serre-Oesterl\'e measure $\op{vol}_{so}$ on the set
    $$(\lie\oplus\lieh)_\pm(O_v).$$
(See \cite{O}.)

\begin{cor}
Let $\lie$ be symplectic or orthogonal.  There exist virtual Chow
motives $M(\sigma)$ over $\Q$, for $\sigma\in\{-1,1\}$, such
that for all number fields $K$ and almost all places $v$ of $K$ we
have
    $$\op{vol}_{so}((\lie\oplus\lieh)(O_v)\cap (\lie\oplus\lieh)_\sigma (F))
    = \op{trace}\,\Frob_v M(\sigma)_K.$$
\end{cor}

\begin{rem}
The trace of Frobenius is to be interpreted as the alternating
trace on the $\ell$-adic cohomology of the motive.
\end{rem}

\begin{proof}
Let $M$ be the virtual Chow motive attached to the virtual set
$$(\lie\oplus\lieh)_\sigma$$ in \cite[3.4]{DL}.   The corollary
is now the comparison theorem of Denef and Loeser, which states
that the Serre-Oesterl\'e measure of definable sets is given by
the trace of Frobenius against the corresponding motive
\cite[8.3.1]{DL}.
\end{proof}

\subsection{Motives attached to unitary Lie algebras}

In the case of the unitary Lie algebra, we obtain a similar
statement, but we must work with the first order language of rings
with involution.  For simplicity, we will take the $p$-adic
extensions $E/F$ defining the unitary Lie algebras to be
unramified.

We must confront the fact that the unitary Lie algebras are not
definable in the first order language of rings.  If $E/F$ is an
unramified quadratic extension of $p$-adic fields of
characterisitic zero, we may fix $\epsilon$ in $E$ such that Galois
conjugation in $E/F$ negates $\epsilon$:
    $$\bar \epsilon = - \epsilon.$$
That is, $\epsilon$ is pure imaginary.   With $\epsilon$ in hand, we may
identify $E$ (basis $1$, $\epsilon$) with the vector space $F^2$ (basis
$(1,0)$, $(0,1)$). Addition and multiplication are replaced with
addition and multiplication expressed in components.  Any formula
in $E$ that involves addition, multiplication, and conjugation can
be replaced with a formula in twice the number of variables
involving component-wise addition, multiplication expressed in
components, and negation on the second factor of each pair $(x,y)$
in $F^2$ representing $x+ y \epsilon$ in $E$.

A similar approach allows us to replace a formula in the
first-order theory of rings with involution with a formula in
Pas's language.  A formula in $n$ free variables
$(z_1,\ldots,z_n)$ becomes a formula in $2n+1$ variables:
        $$z_i \mapsto x_i +  y_i\,\epsilon$$
By equating real and imaginary parts, we may eliminate the variable $\epsilon$
from the equations, leaving only
formulas involving $u = \epsilon^2$.  The free variables of the resulting
formula are $u$,  $x_i$, and $y_i$.

It is not possible to constrain $u$ within the first order
language of rings to be a particular element that is not a square.
($-1$, for example, is a square in some rings and not in others,
and we do not want to restrict the $p$-adic domains by assuming
that it is a square.) But within Pas's language, we can constrain
$u$ to lie within a given definable set of nonsquares. Thus, we
define $\phi(u)$ to be the formula
with quantifier $\forall \xi$ of the residue field sort:
        $$\forall \xi (\xi^2\ne \op{ac}(u)) \wedge (\op{ord}(u)=0).$$
This constrains $u$ within a set that yields isomorphic
unramified field extensions
$F(\sqrt{u^F})$ for all interpretations $u^F$ of $u$ in $F$.
The
involution-free formula in $2n+1$ variables defines a
$1$-parameter family of structures. The unitary Lie
algebra, for example, is replaced by a $1$-parameter family of
isomorphic Lie algebras.  Although each algebra is not
individually definable, the $1$-parameter family is.

To the formula in $2n+1$ variables in Pas's language (conjoined
with the constraint $\phi(u)$), we attach a virtual Chow motive as
before over $\Q$.  The trace of Frobenius against this
motives computes the Serre-Oesterl\'e measure within the
$1$-parameter family of isomorphic objects.

\section{$p$-adic theory}

Let $F$ be a $p$-adic field; following Waldspurger, we
assume that the residue field
characteristic $p$ is sufficiently large ($p\ge 3\dim(\lie)+1$).

In contrast to the previous section, here we  consider the Lie
algebras over the field $F$, not their virtual counterparts. We
will use the notation of set theory in its traditional meaning
here. We still think of the classical Lie algebras as defined by
means of the same matrices $J$ as in \ref{liealgebras}; however,
in this section, elements are to be understood as actual matrices with
entries in the $p$-adic field, not formal symbols.

\subsection{Parametrization of  regular semisimple orbits}\label{orbits}
We will need a parametrization of regular semisimple orbits in the
classical Lie algebras. Here we quote it in detail from
Waldspurger, \cite[I.7]{W}.

Denote by $(V,q_V)$ the vector space on which $\lie$ acts by
endomorphisms, with the
quadratic form preserved by $\lie$.
First, consider  the case $\lie=\syp(2r, F)$ or $\lie=\so(2r+1,
F)$. Let $X$ be a regular semisimple element in $\lie(F)$. Then
the orbit of $X$ corresponds to the following data:
\begin{itemize}
\item
    a finite set $I$
\item
    a finite extension $F_i^{\#}$ of $F$ for each $i\in I$,
    \item
    a
    2-dimensional $F_i^{\#}$-algebra $F_i$,
    \item an element $a_i\in
    F_i^{\times}$,
\item
    a collection of elements $c_i\in F_i^{\times}$, $i\in I$,
\end{itemize}
subject to the following conditions:
\begin{itemize}
\item
    For $i,j\in I$, $i\neq j$,
    there is no $F$-linear isomorphism
    between $F_i$ and $F_j$ taking $a_i$ to $a_j$
\item
    $a_i$ generates $F_i$ over $F$
\item
        For all $i\in I$, denote by $\tau_i$ the unique nontrivial
        automorphism  of $F_i/F_i^{\#}$. Then $\tau_i(a_i)=-a_i$.
\item
  $\sum_{i\in I}[F_i:F]=\dim V=2r$
        in the case $\lie=\syp(2r)$,
  $\sum_{i\in I}[F_i:F]=\dim V-1=2r$
        in the case $\lie=\so(2r+1)$
\item
  $\tau_i(c_i)=-c_i$
      in the symplectic case;
  $\tau_i(c_i)=c_i$
        in the orthogonal case
\item
    In the symplectic case: set
        $W=\bigoplus_{i\in I}F_i$
    and let $X_W$ be the element of $\op{end}\,W$ defined by
    $X_W(\sum_{i\in I}w_i)=\sum_{i\in I}a_iw_i$.
    Then $(V, q_V)$ is isomorphic to the space $W$
    endowed with the form
$$
   q_W\left(\sum_{i\in I}w_i, \sum_{i\in I}w_i'\right)= \sum_{i\in
   I}\op{trace}_{F_i/F}(\tau_i(w_i)w_i'c_i).
$$
    The isomorphism between
    $V$ and $W$ allows to identify $\op{end}\,W$ with $\op{end}\,V$, and $X_W$
    is identified with some $X\in \lie(F)$.
    (Note that our constants $c_i$ differ from
    the ones in \cite{W} by a factor of $[F_i:F]^{-1}$.)

        In the odd orthogonal case,
    define $(W, q_W)$ as in the previous
    case. There is an additional requirement that there exists a
    one-dimensional orthogonal space $(W_0,q_0)$ over $F$, such that
    $(W_0\oplus W,q_0\oplus q_W)$ is isomorphic to $(V, q_V)$.
    The action of the element $X_W$ on $W$ is defined as above, and
    $X_W$ acts by zero on $W_0$.
\item
             In the orthogonal case,
    when
    the space $(W_0,q_0)$ exists, its class is determined uniquely.
\end{itemize}

A different choice of the isomorphism between $(V, q_V)$ and
$(W,q_W)$ would identify $X_W$ with an element in the same $F$-conjugacy
class.
The orbit ${\mathcal O}(X)$ is well-defined and is
denoted $\orb(I,(a_i),(c_i))$. The correspondence between the
orbits and the data is one-to-one if we identify the triples
$(I,(a_i),(c_i))$ and $(I',(a_i'),(c_i'))$ subject to the
following conditions.
\begin{itemize}
\item
    There is a bijection $\phi:I\to I'$
\item
    For all $i\in I$ there is an $F$-linear isomorphism
    $\sigma_i\colon F_{\phi(i)}'\to F_i$ such that
    $\sigma_i(a_{\phi(i)}')=a_i$
\item
    For all $i\in I$, denote by
    $\sgn_{F_i/{F_i^{\#}}}$ the quadratic character of
    $F_i^{\#}$ associated with the algebra $F_i$. Then
    $\sgn_{F_i/F_i^{\#}}(c_i\sigma_i({c_{\phi(i)}'}^{-1}))=1$.
    (Notice
    that by definition of the $c_i$ and $c_i'$, the product
    $c_i\sigma_i({c_{\phi(i)}'}^{-1})$ is stable under $\tau_i$, and
    therefore lies in $F_i^{\#}$).
\end{itemize}

The {\it stable} orbit of $X$ does not depend on the constants $c_i$.

In the case $\lie=\so(2r)$, the parametrization of the orbits
needs to be modified as follows. The data $(I,(a_i),(c_i))$ are
defined in the same way as in the odd orthogonal case. However,
the correspondence between the data $(I, (a_i),
(c_i))$ and the orbits is no longer one-to-one. First, the
construction gives only the orbits that do not have the
eigenvalue $0$. Second, depending on the choice of the
isomorphism between $(W, q_W)$ and $(V, q_V)$, the element $X_W$
maps into one of the two distinct orbits, which remain distinct
even over the algebraic closure of $F$. They will be  denoted
by $\orb^+(I,(a_i),(c_i))$ and $\orb^-(I,(a_i),(c_i))$. Their
union is denoted by $\orb(I, (a_i),(c_i))$.

Suppose $X\in \so(2r, F)$, $X\in\orb(I, (a_i),(c_i))$ for some
data $(I, (a_i),(c_i))$. It is possible to tell whether
$X\in\orb^+$ or $X\in\orb^-$, using the Pfaffian of the matrix
$JX$. Assume that $A$ is skew-symmetric.
Recall the properties of Pfaffian that we need
\cite[B.2.6]{G}:
\begin{itemize}
\item
    $\op{pf}(A)^2=\det A$.
\item
    If $g\in GL(2r,F)$, then
    $\op{pf}({}^tgAg)=\det(g)\op{pf}(A)$.
\item
    Let $A$ and $B$ be two skew-symmetric matrices of sizes $2a$ and
    $2b$, respectively. Let $A\oplus B$ be the block-diagonal matrix
    with the diagonal blocks $A$ and $B$.
    Then $\op{pf}(A\oplus B)=\op{pf}(A)\op{pf}(B)$.
\end{itemize}

We will also need the following embedding $\psi$ of
$\so(2a,F)\oplus\so(2b,F)$ into $\so(2(a+b),F)$: for
$A\in\so(2a,F)$, $B\in\so(2b, F)$ let $\psi(A,B)$ be the block matrix
\begin{equation*}
 \left(
 \begin{matrix}
    B_1            & 0 & B_2\\
    0 & A      & 0 \\
    B_3            & 0 & B_4
 \end{matrix}
\right), \quad {\text{where}} \quad
 B=
  \left(
  \begin{matrix}
    B_1 & B_2\\
    B_3 & B_4
  \end{matrix}
 \right)
\end{equation*}

\begin{lem}\label{pfaff}
1. The function $\theta(X)=\op{pf}(JX)$ is constant on
each one of the sets $\orb^+(I, (a_i),(c_i))$ and
$\orb^-(I,(a_i),(c_i))$ in $\so(2r, F)$, and takes distinct values
on these two sets.

2. Let $\psi$ be the embedding of $\so(2a,F)\oplus\so(2b,F)$ into
$\so(2(a+b), F)$. Then $\theta(\psi(A,B))=\theta(A)\theta(B)(-1)^{a b}$.
\end{lem}

{\bf Proof.} {\bf 1.} Suppose $X, X'\in\orb(I, (a_i),(c_i))$, and
suppose $X=g^{-1}X'g$. Then $JX=Jg^{-1}X'g={}^tgJX'g$, so
$\op{pf}(JX)= \det(g)\op{pf}(JX')$. Hence, $\op{pf}(JX)=\op{pf}(JX')$ if and only
if $\det g=1$. That is, $X$ and $X'$ are conjugate within
the special orthogonal group.

{\bf 2.} In order to use the properties of the Pfaffian, we need
to bring the matrix $\psi(A,B)$ to block-diagonal form. This is
done by conjugating in $GL_{2a+2b}(F)$ by the
permutation matrix
\begin{equation*}
 \sigma=\left(
\begin{matrix}
 0        & \text{Id}_{2a\times 2a} & 0\\
 \text{Id}_{b\times b} &  0          & 0 \\
 0         & 0          & \text{Id}_{b\times b}
\end{matrix}
\right).
\end{equation*}
Note $J_{2a+2b} = \psi((-1)^b J_{2a},J_{2b})$ by Definition \ref{Jeven}.
Then
$$\begin{array}{lll}
\theta(\psi(A,B)) &= \op{pf}(J_{2a+2b}\psi(A,B))=
\op{pf}(\psi((-1)^bJ_{2a}A,J_{2b}B))\\
    &=\det \sigma\,\op{pf}((-1)^bJ_{2a}A\oplus J_{2b}B)
=\theta(A)\theta(B) (-1)^{ab},
\end{array}$$
since $\det \sigma=1$. \qed

The remaining case is the case of the unitary group,
$\lie=\un(n)$. The orbits in $\un(n)$ are parametrized
similarly, with the data $(I,(a_i),(c_i))$. However,
$F_i$ is obtained by a different construction.
$F_i^{\#}$ is
an extension of the base field $F$ as before, but  $F_i$ now
equals the algebra  $F_i^{\#}\otimes_F E$, where $E$ is the quadratic
extension of $F$, used to define
the given unitary Lie algebra.
In the unitary case, our form differs from that in Waldspurger by
a constant $\eta\in E$ that he chooses.
Our constant $c_i$ also differs from his
by a factor of $\eta$.

As before, in the unitary case  the parameters $a_i$ and $c_i$ lie
in $F_i^{\times}$; they are subject to the following conditions:
\begin{itemize}
\item
    For any $i\in I$, $a_i$ generates $F_i$ over $E$
\item
    There is no $E$-linear isomorphism of $F_i$ onto $F_j$ taking
    $a_i$ to $a_j$ for $i\neq j$
\item
    Let $x\mapsto \bar{x}$ denote the unique nontrivial element
    of $\op{Gal}(E/F)$;
    and for each $i\in I$ let $\tau_i$ be the automorphism
    ${\op{id}}\otimes\ \bar{(\cdot)} $.
    Then $\tau_i(a_i)=-a_i$, $\tau_i(c_i)=c_i$, $i\in I$.
\item
    $\sum_{i\in I}[F_i:E]=r$.
\end{itemize}

As in the other cases, set $W=\bigoplus_{i\in I}F_i$. It is a
vector space over $E$. The hermitian form $q_W$ on $W$ is defined
by
$$
  q_W\bigl(\sum_{i\in I}w_i,\sum_{i\in I}w_i'\bigr)
= \sum_{i\in I} {\op{trace}}_{F_i/E}(\tau_i(w_i)w_i'c_i).
$$
Note that as before, our $c_i$ differ from Waldspurger's by a
factor of $[F_i:E]$. The rest of the construction remains the
same; it gives all orbits, and the correspondence is one-to-one
with the same definition of data $(I,(a_i),(c_i))$ and
$(I',(a_i'),(c_i'))$ as in the symplectic case, except to
change `$F$-linear isomorphism' to `$E$-linear isomorphism'.

\subsection{The function $\sign$}\label{thefactor}
Recall the Waldspurger's description of the $\sign$ function in
detail. Here we are again quoting from \cite[X.8]{W}. Let
$P_X\in F[\lambda]$ be the characteristic polynomial of
$X$, as before. Let $\sign$ be the function on
$\lie\oplus{\mathfrak h}$ from the definition of the transfer
factor.

The function $\sign(X,Y,Z)$ takes the value $0$ unless $X\in \lie$
and $(Y,Z)\in {\mathfrak h}$ are regular semisimple, and the
stable conjugacy classes of $X$ and $Y$ correspond. This means
that the condition of Definition
\ref{pmset} is fulfilled (cf.  Lemma \ref{pfaff}).

Suppose that $X$ and $(Y,Z)$ are regular semisimple elements of
$\lie(F)$ and ${\mathfrak h}(F)$, respectively, and their stable
conjugacy classes correspond. Suppose that
    $Y\in\orb(I_1,(a_{i,1}),(c_{i,1}))$,
    $Z\in \orb(I_2, (a_{i,2}),(c_{i,2}))$
(Section \ref{orbits}). Let $I=I_1\cup I_2$, and let $(a_i)$ be
the join of the lists $(a_{i,1})$ and $(a_{i,2})$. Since the
stable conjugacy classes of $(Y,Z)$ and $X$ correspond, there
exists a family $(c_i)$, such that $X\in\orb(I,(a_i),(c_i))$.
Define constants $C_i$ (upper case) by
$$C_i=\begin{cases}
    c_i^{-1}a_i^{-1}P_X'(a_i) & \text{if $\lie$ is even orthogonal}\\
    c_i^{-1}P_X'(a_i)& \text{otherwise.}
    \end{cases}.$$
The elements $C_i$ lie in $F_i^{\#}$ for all $i\in I$.
\begin{thm}
(Waldspurger, \cite{W}, Section X.8)
$$
  \sign(X,Y,Z)=\prod_{i\in I_2^{\ast}}\sgn_{F_i/F_i^{\#}}(C_i),
$$
where $I_2^{\ast}$ is the set of all indices $i\in I_2$ such that
$F_i/F_i^{\#}$ is a
field extension.
\end{thm}

\section{Proof of the main theorem}

\subsection{A few $p$-adic lemmas}\label{embedding}

As above, let $F$ be a $p$-adic field, $I$ -- a finite set,
$F_i^{\#}$, ${i\in I}$ -- a collection of finite extensions of
$F$, $F_i$ -- a 2-dimentional $F_i^{\#}$-algebra for every $i\in
I$. Let $\tau_i$ be the unique nontrivial automorphism of
$F_i/F_i^{\#}$, and let $\sgn_{F_i/{F_i^{\#}}}$ be the quadratic
character of $F_i^{\#}$ associated with the algebra $F_i$. The
indexing set $I$ is a union of two sets $I^{\ast}$ and $I^{0}$,
where $F_i$ is a field extension of $F$ for $i\in I^{\ast}$, and
$F_i$ is a direct sum of two copies of $F_i^{\#}$ for $i\in
I^{0}$. For $i\in I^{\ast}$, $\tau_i$ is the nontrivial
Galois automorphism of $F_i/F_i^{\#}$; for $i\in I^{0}$, $\tau_i$
exchanges the two copies of $F_i^{\#}$ in $F_i$. We observe that
 the character $\sgn_{F_i/{F_i^{\#}}}$ is nontrivial if and only
if $i\in I^{\ast}$.

Let $W=\bigoplus_{i\in I}F_i$. Let $\phi$ be the isomorphism from Section
\ref{orbits} between $(W,q_W)$ (or, in the odd orthogonal case,
$(W\oplus W_0, q_W\oplus q_{W_0})$) and $(V,q_V)$. Let
$\phi_{\ast}$ be the isomorphism
    $\phi_{\ast}\colon \op{end}\,W\to\op{end}\,V$
induced by $\phi$. The map $\phi_{\ast}$ induces an isomorphism from
centralizer of $X_W$ in $\op{end}\,W$ onto
$C(X)$.

Let $L_i\colon F_i\to \op{end}(F_i)$ be the linear map that takes
$w\in F_i$ to the operator that acts by
multiplication by $w$ on the left on $F_i$. (Generally,
 $\op{end}$ should be understood as $\op{end}_F$; except, in the
unitary case, where it should be understood as
 $\op{end}_E$.) Let
$L\colon W\to\op{end}\,W$ be the direct sum of the maps $L_i$.

\begin{defn}
Let $X$ be a  regular semisimple element in a classical Lie
algebra $\lie$. Define the involution $\tau$ on $C(X)$ by
    $\tau(g)=(J^{-1}){}^tgJ$
with the appropriate matrix $J$ in the symplectic and orthogonal
cases, and by
    $\tau(g)=(J^{-1}){}^t\bar{g}J$ in the unitary case.
\end{defn}

\begin{lem}\label{tau}
Let $\tau_i$ be the involutions on $F_i$ from Section
\ref{orbits}. Let $\tau_W$ be the involution on $W$ that acts by
$\tau_i$ on each $F_i$, $i\in I$. For all $a\in W$,
$$
 \phi_{\ast}\circ L(\tau_W(a))=\tau(\phi_{\ast}\circ L(a)).
$$
\end{lem}\nd
{\bf Proof.}
First, observe that $\dim W=\dim C(X)$. It follows that
$\phi_{\ast}\circ L$ is an isomorphism of $F$-algebras
($E$-algebras in the unitary case). Consider the two involutions
on $C(X)$: $\tau$ and
    $\tau_{\ast}=(\phi_{\ast}\circ L)\circ {\tau_W}
        \circ (\phi_{\ast}\circ L)^{-1}$.
We show that they coincide.

It will be more convenient to consider $-\tau$ and $-\tau_{\ast}$. It
suffices to show that the involutions have the same set of fixed points. The
set of fixed points of $-\tau$ is precisely $C_{\lie}(X)$. Let $Y\in C(X)$
be a fixed point of $-\tau_{\ast}$, and let
    $a=(\phi_{\ast}\circ L)^{-1}(Y)$.
Then $\tau_W(a)=-a$, since $Y$ is a fixed point  of $-\tau_{\ast}$. By
definition of the form $q_W$, the condition $\tau_W(a)=-a$ is
equivalent to the condition
$$
q_W(L(a)w',w'')+q_W(w',L(a)w'')=0\quad\text{for all}\quad
w',w''\in W.
$$
Since $\phi(q_W)=q_V$, the latter condition is equivalent to
$\tau(Y)=-Y$;
that is, to the condition $Y\in C_{\lie}(X)$. \qed

We will need the following observation.
\begin{rem}\label{c}
Let $w=(c_i)_{i\in I}\in W$ be part of the data from Section
\ref{orbits}. Then $c=\phi_{\ast}\circ L(w)$ is an invertible
element in $C(X)$, possessing the following property.

{\it There exists a basis $e_1,\dots,e_r$ of $C(X)$ such that for
all $x=(x_1,\dots,x_r)$ and $x'=(x_1',\dots,x_r')$}
$$\op{trace}\left(\tau\bigl(\sum_{i=1}^rx_ie_i\bigr)
\bigl(\sum_{i=1}^r x_i'e_i\bigr)c\right)={}^tx Jx'$$
(or ${}^t\bar x J x$, if unitary).
Indeed, the existence of the basis is the same condition as the
existence of the isomorphism $\phi$. The matching of the quadratic
forms follows from the definition of $\phi$, the previous lemma,
and the fact that trace of an element of a field is equal to the
trace of the endomorphism defined by left multiplication by that element. The
invertibility of $c$ follows from the nondegeneracy of the
quadratic form defined above.
\end{rem}

\begin{lem}\label{projsigns}
Let $I$ be part of the data for $C(X)$ as above, and let $i\in I$.
Let $z=z_i\in F_i^{\#}$ be an arbitrary element. Then the
following conditions are equivalent:
\begin{enumerate}
\item\label{sg}
    $\sgn_{F_i/F_i^{\#}}(z_i)=1$
\item\label{nm}
    Let $w$ be an arbitrary element of $W$ such that
    $w_i=z_i$.
    Let $P$ be the projector from $C(X)$ onto
    $\phi_{\ast}\circ L(F_i)$.
    Then there exists $X_1\in C(X)$ such that
    $P X_1\tau(X_1)=P (\phi_{\ast}\circ L) (w)$.
\end{enumerate}
\end{lem}\nd
{\bf Proof.} If $F_i/F_i^{\#}$ is a field extension
and if $z_i\in F_i^{\#}$, then $\sgn_{F_i/F_i^{\#}}(z_i)=1$ if and
only if $z_i$ is a norm of an element of $F_i$. The condition that
an element $z_i$ is a norm can be written as $\exi y\in F_i\colon
z_i=\tau_i(y)y$. In the case when $F_i$ is an algebra, both
conditions $\sgn_{F_i/F_i^{\#}}(z_i)=1$ and $z_i=\tau_i(y)y$ for
some $y\in F_i$ hold for all $z_i\in F_i$.

Suppose the condition (\ref{sg}) holds. Then let $y_i$ be the
element of $F_i$ such that $\tau_i(y_i)y_i=z_i$. Let $y$ be the
element of $W$ whose $j$-th component is $y_i$ if $j=i$, and $0$
if $j\in I\setminus \{i\}$. Let $X_1=\phi_{\ast}\circ L(y)$. Then,
by Lemma \ref{tau}, $\tau(X_1)=\phi_{\ast}\circ L(\tau_W(y))$.
Hence, $\tau(X_1)X_1=\phi_{\ast}\circ L(\tau_W(y)y)$. The
observation that $P (\phi_{\ast}\circ L)(\tau_W(y)y)= P
(\phi_{\ast}\circ L)(w)$ since $\tau_W(y)y$ and $w$ have the same
$i$-th component, completes the proof of the implication
(\ref{sg})$\Rightarrow$(\ref{nm}). The converse is proved by
reversing all the steps. \qed

\subsection{Proof of the Main Theorem}\label{proof}
Let $\lie\oplus\lieh$ be a virtual Lie algebra as in Section
\ref{liealgebras}. Let $F$ be a $p$-adic field of sufficiently
large residue field characteristic.
 Let
    $\lie(F)\oplus\lieh(F)$
be the Lie algebra over $F$ attached to the virtual Lie algebra
$\lie\oplus\lieh$. The transfer factor $\Delta(X,Y,Z)$ is defined
on
    $\lie(F)\oplus\lieh(F)$.

First, consider the zero set of $\Delta$, that is, $\sigma=0$.
It is immediately clear by comparison of the definition of the
transfer factor in Section \ref{thefactor} with Definiton \ref{pmset}
that $\Delta(X,Y,Z)\neq 0$ if and only if  $(X,Y,Z)$ satisfies all
the conditions of Definition \ref{pmset}, that is,
    $(X,Y,Z)\notin (\lie\oplus\lieh)_0(F)$.

Now let $\sigma =1$. Suppose that
    $(X,Y,Z)\in (\lie\oplus\lieh)_{\pm}(F)$.
We need to check that  $\sign(X,Y,Z)=1$ if and only if
the formula from Section \ref{1-set} evaluated at $(X,Y,Z)$
is true in the interpretation provided by $F$.
The virtual centralizer ${\bf C}(X)$ becomes the
actual centralizer of the element $X$ in $\gl(n,F)$. We keep the
notation of Section \ref{orbits}.

Let $Y\in\orb(I_1,(a_{i,1}),(c_{i,1}))$, $Z\in
\orb(I_2,(a_{i,2}),(c_{i,2}))$. Recall from Section
\ref{thefactor} that then there exists an element $(c_i)\in W$
such that $X\in\orb(I,(a_i),(c_i))$, where $I=I_1\cup I_2$, and
$(a_i)$ is the direct sum of $(a_{i,1})$ and $(a_{i,2})$
as
elements of the vector space $W$.

By Lemma \ref{tau} and Remark \ref{c}, the existence of $(c_i)_{i
\in I}$ such that $X\in\orb(I,(a_i),(c_i))$ and the $c_i$ satisfy
all the requirements of Section \ref{orbits}, is equivalent to the
conjunction of the following statements:
\begin{itemize}
\item The  elements $c_i$ satisfy $\tau_i(c_i)=\chi c_i$, where
    $\chi=-1$ if $\lie$ is symplectic or even unitary,
    and $\chi =1$ otherwise.
      Let $c=\phi_{\ast}\circ L((c_i)_{i\in I})$.
      By Lemma \ref{tau}, this condition is equivalent to
     $\tau(c)=\chi c$.
\item
  There exists $w=(c_i)\in W$
  and a basis $e_1,\dots, e_n$ of $C(X)$, satisfying the conditions
  of Remark \ref{c}.
  This condition can be written in the terms of logic
  as
    $$\exists c \in {\bf C}(X)\quad \op{trace-form}(X,c)$$
\end{itemize}

Now we need to check that the condition of Waldspurger
$$
\prod_{i\in I_2^{\ast}}\sgn_{F_i/F_i^{\#}}(C_i)=1
$$
is equivalent to our condition
$\op{even-parity}(P_{Z,0},\phi(\cdot,X,c'))$.

In the unitary case,
observe that we may take the product over all of $I_2$
instead of $I_2^{\ast}$, since the $\sgn$ character is trivial if
$F_i/F_i^{\#}$ is not a field extension. Let $X_W=(a_i)_{i\in
I}\in W$.
In the non-unitary case, an irreducible factor of $P_{Z,0}$ is an
even polynomial iff the corresponding algebra $F_i$ is a field.
Thus, the antecedent
    $$\op{even-poly}(f) \implies \cdots $$
correctly distinguishes between $I_2^{\ast}$ and $I_2$.

Since $\phi_{\ast}\circ L$ is an isomorphism
of $F$-algebras ($E$-algebras in the unitary case), we have
$(\phi_{\ast}\circ L)(p(X_W))=p(\phi_{\ast}\circ L(X_W))$ for any
polynomial $p$ with coefficients in $F$.
We will use this observation with $p=P_X'$.

Let $M=\phi_{\ast}\circ L((C_i)_{i\in I})$
(it is well defined since $C_i\in F_i^{\#}$).
In the symplectic, unitary, and odd orthogonal
cases,
$$
M=\phi_{\ast}\circ L((c_i^{-1}P'_X(a_i))_{i\in I})
= c'P_X'(X).
$$
The last equality holds by the definition of $c$, $c'$, and
$X=\phi_{\ast}\circ L((a_i)_{i\in I})$. In the even orthogonal
case, the right hand side of the above formula will have an extra
factor of $X^{-1}$.

Suppose that $\prod_{i\in I_2}\sgn(C_i)=1$.
Then
$\sgn(C_i)=-1$ for an even number of
indices $i\in I_2$. Let $P_{Z,0}(\lambda)=f_1\dots f_m$ (here $m$
equals the cardinality of $I_2$) be the factorization of the
characteristic polynomial of $Z$ into irreducibles. By Lemma
\ref{projsigns}, for each $i=1,\dots ,m$, the condition
$\sgn_{F_i/F_i^{\#}}(C_i)=1$ is equivalent to the existence of a
matrix $X_i$ such that
$\Pi(X,f_i,\tilde f_i)\tau(X_i)X_i=\Pi(X,f_i,\tilde f_i)M$.
The calculation of $M$ shows that
$\sgn_{F_i/F_i^{\#}}(C_i)=1$ if and only if the condition
$\op{norm}(X,f,\ldots)$ holds, where $c'c=1$. By definition of the
condition $\op{even-parity}(P_{Z,0},\phi(\ldots))$,
it holds if and only if the
number of $\op{norm}(X,f,\ldots)$ that fail is even,
that is, if and only
if $\Delta(X,Y,Z)$ is positive.
(In the unitary case, we must take the factorizations to be of $P_{X\epsilon}$
and $P_{Z\epsilon}$.)
This completes the proof.

\end{document}